\documentclass[graybox]{svmult}

\usepackage{amsmath,amssymb}
\usepackage{mathptmx}        
\usepackage{helvet}          
\usepackage{courier}         
\usepackage{type1cm}        

\usepackage{graphicx}        
\usepackage{multicol,multirow} 
\usepackage[bottom]{footmisc}
\usepackage{latexsym,bm}        
\usepackage{algorithm,algorithmic,tikz}
\usepackage[labelformat=simple]{subfig}

\graphicspath{ {./eps/} }

\newcommand{\Cross}{\mathbin{\tikz [x=1.4ex,y=1.4ex,line width=.2ex] \draw (0,0) -- (1,1) (0,1) -- (1,0);}}%

\usepackage{lineno}

\begin{document}

\title*{Energy stable model order reduction for the  Allen-Cahn equation}
\titlerunning{Energy stable MOR for Allen-Cahn equation}

\author{Murat Uzunca \and  B\"{u}lent Karas\"{o}zen}
\authorrunning{Uzunca and Karas\"{o}zen}

\institute{M.~Uzunca \at Department of Industrial Engineering, University of Turkish Aeronautical Association \& Institute of Applied Mathematics, Middle East Technical University, Ankara, Turkey, \\ \email{muzunca@thk.edu.tr}
\and  B.~Karas\"{o}zen \at Institute of Applied Mathematics \& Department of Mathematics, Middle East Technical University, Ankara, Turkey, \\ \email{bulent@metu.edu.tr}
}

\maketitle

\abstract{The Allen-Cahn equation is a gradient system, where the free-energy functional decreases monotonically in time. We develop an energy stable reduced order model (ROM) for a gradient system, which inherits the energy decreasing property of the full order model (FOM).  For the space discretization we apply a discontinuous Galerkin (dG) method and for time discretization the energy stable average vector field (AVF) method. We construct ROMs  with proper orthogonal decomposition (POD)-greedy adaptive sampling of the snapshots in time and evaluating the nonlinear function with greedy discrete empirical interpolation method (DEIM). The computational efficiency and accuracy of the reduced solutions are  demonstrated numerically for the parametrized Allen-Cahn equation with Neumann and periodic boundary conditions.
\keywords{Gradient systems, energy stability, discontinuous Galerkin, proper orthogonal decomposition, empirical interpolation, adaptive sampling.}
}

\section{Introduction}
\label{intro}

The Allen-Cahn equation  \cite{allen79amt}
\begin{linenomath*}
\begin{equation}\label{allencahn}
u_t=\epsilon\Delta u -f(u), \quad (x,t) \in \Omega \times (0,T],
\end{equation}
\end{linenomath*}
on a bounded region $\Omega\subset\mathbb{R}^d (d=1,2)$, is a gradient system in the $L_2$ norm:
\begin{linenomath*}
\begin{equation} \label{acgradient}
u_t = - \frac{\delta \mathcal{E}(u)}{\delta u}.
\end{equation}
\end{linenomath*}
Equation \eqref{acgradient} is characterized by the minimization of the Ginzburg--Landau energy functional
\begin{linenomath*}
\begin{equation*}
\mathcal{E}(u)=\int_{\Omega} \left( \frac{\epsilon}{2}|\nabla u|^2 + F(u)  \right)dx,
\end{equation*}
\end{linenomath*}
with a potential functional $F(u)$. The main characteristic of a gradient system is the energy decreasing property:
\begin{linenomath*}
\begin{equation}\label{energydecrrease}
\mathcal{E}(u(t_{n})) < \mathcal{E}(u(t_{m})), \quad \forall t_{n} > t_{m}.
\end{equation}
\end{linenomath*}

The Allen-Cahn equation \eqref{allencahn} was originally introduced to describe the phase of a binary mixture. Nowadays it is used as a model for interface problems in material science, fluid dynamics, image analysis, mean curvature flow, and pattern formation. In \eqref{allencahn} the unknown $u$ denotes the concentration of the one of the mixture. The parameter $\epsilon$ is  related to the interfacial  width, capturing the dominating effect of reaction kinetics and stays for effective diffusivity. The non-linear term $f(u)$ in \eqref{allencahn} is given by $f(u)=F'(u)$.  Depending on the choice of the potential functional $F(u)$, i.e. the non-linear function $f(u)$, different types of gradient systems occur. The most common potential functions for the  Allen--Cahn equation are  the convex quartic double-well potential \cite{feng13nsi} and the non-convex logarithmic potential \cite{barrett99fea}, given respectively by:
\begin{linenomath*}
\begin{subequations}
\begin{align}
F(u) &= (u^2-1)^2/4, \label{quartic} \\
F(u) &= ( \theta[(1+u)\ln(1+u) + (1-u)\ln(1-u)] - \theta_cu^2 )/2, \label{logarithmic}
\end{align}
\end{subequations}
\end{linenomath*}
where $\theta_c$ in \eqref{logarithmic} is the transition temperature. For temperature $\theta$ close to $\theta_c$, the logarithmic potential is usually approximated by the convex quartic double-well potential.  In case of the quartic double-well potential, $f(u)=u^3-u$  represents the bi-stable nonlinearity. For the logarithmic potential, it takes the form $f(u)=(\theta/2)\ln\left((1+u)/(1-u)\right) - \theta_c u$.

The main characteristic of the Allen-Cahn equation  \eqref{allencahn} is the rapid formation of the transient layers and exponentially slow formation of the terminal layers for very small values of $\epsilon$. This is known as metastability phenomena, characterized by the relative flatness of solutions, where the stable or unstable fixed points coalesce or vanish over long time. These make the numerical computation of the Allen-Cahn equation \eqref{allencahn} challenging for very small values of $\epsilon$. In the literature for discretization of \eqref{allencahn} in space, the well-known finite-differences, spectral elements \cite{christlieb14has}, continuous finite elements \cite{liu13ssi} and local discontinuous Galerkin (LDG) method \cite{Guo16} are used. Several energy stable integrators are developed to preserve the energy decreasing property of the Allen-Cahn equation. For small values of the diffusion parameter $\epsilon$, semi-discretization in space leads to stiff systems. Therefore it is important to design efficient and accurate numerical schemes that are energy stable  and robust for small $\epsilon$. Because the explicit methods are not suitable for stiff systems, several energy stable implicit-explicit methods based on the convex splitting of the non-linear term are developed \cite{shen10nac,choi09ugs,feng13scn}.

In this work, we use the symmetric interior penalty Galerkin (SIPG) finite elements for space discretization  \cite{arnold82ipf,riviere08dgm} and the energy stable average vector field (AVF) method  \cite{celledoni12ped,hairer10epv} for time discretization. The SIPG approximation enables to capture the sharp gradients or singularities locally.  On the other hand, the AVF method is the only second order implicit energy stable method for a gradient system. Because the computation of the patterns for small  values of $\epsilon$  is time-consuming, we consider reduced order modeling which inherits the essential dynamics like the energy decreasing property of the Allen-Cahn equation. In the literature, there are only two papers dealing with the reduced order modeling for Allen-Cahn equation. Using finite difference discretization in space and convex splitting in time, an energy stable reduced order model is derived in \cite{Song16}. In \cite{Kalashnikova12}, a non-linear POD/Galerkin reduced order model is applied for efficient computation of the metastable states. An early application of POD to the optimal control of phase field models in material sciences dates back to 2001 \cite{Volkweinn01}. It was shown that for the optimal control of two coupled non-linear PDEs, the solutions of the POD reduced model have nearly the same accuracy as the finite element FOM solutions, whereas the computing time is reduced enormously. Here we apply the greedy proper orthogonal decomposition (PODG) method for the parametrized non-linear parabolic PDEs \cite{Drohmann12,Grepl13}, where the reduced basis functions are formed iteratively by a greedy algorithm for the parameter values such that a POD mode from the matrix of projection error for the parameter value with the largest error is captured and used to enlarge the reduced space.  In order to reduce the computational complexity of the function evaluation for the non-linear term in the reduced model, the empirical interpolation method (EIM) \cite{Barrault04,Grepl07} and discrete empirical interpolation method (DEIM)  \cite{chaturantabut10nmr,Wirtz14} are used. The greedy DEIM is included in the adaptive sampling algorithm. We use in the PODG sampling algorithm, the residual-based a posteriori error indicator such that the FOMs are solved only for  selected parameter values. We see that the DEIM reduced system is conditionally energy stable whereas the fully discrete system is unconditionally stable. The performance of the PODG approach is illustrated for the Allen Cahn equation with quartic and logarithmic potential functions for different parameters. We want to remark that the majority of the model order reduction (MOR) techniques for parametrized PDEs are projection based, which use the state-space description of the models by numerical simulation. There exists equation-free MOR methods using system responses such as measurements \cite{Benner15}. The data-driven MOR method in the Loewner framework \cite{Ionita14} was  applied to parametrized systems.

The rest of the paper is organized as follows.  In Section \ref{dgfem}, fully discretization  of the  model problem \eqref{allencahn} is introduced. In Section \ref{sec:rom}, we describe the reduced order modeling together with a POD-greedy sampling algorithm. We present in Section \ref{numerical} numerical results for ROMs of the parametrized Allen-Cahn equation with the parameters $\epsilon$ and $\theta$.

\section{Fully discrete system}
\label{dgfem}

In this section, we describe the full discretization of the parametrized form of the (2D)  Allen-Cahn equation \eqref{allencahn} using the symmetric interior penalty Galerkin (SIPG) in space and the second order energy stable average vector field (AVF) method in time. For a certain parameter $\mu$, we denote the parameter dependence of a solution $u(x,t)$ by $u(\mu):=u(x,t;\mu)$, where the parameter $\mu$ stands for either the diffusivity $\epsilon$, or the temperature $\theta$ in case of logarithmic potential. We also denote the parameter dependence of the non-linear function by $f(u;\mu)$. Then, the variational form of the parametrized Allen-Cahn equation is given as:
\begin{linenomath*}
\begin{align}\label{varformallencahn}
(\partial_{t}u(\mu), \upsilon)_{\Omega} + a(\mu;u(\mu), \upsilon) + (f(u(\mu);\mu ),\upsilon)_{\Omega}&=0, & \forall \upsilon \in H^1(\Omega).
\end{align}
\end{linenomath*}

The Allen-Cahn equation was considered in the literature under Dirichlet, Neumann and periodic boundary conditions. Here, we give the SIPG discretization for the homogeneous Neumann boundary conditions \cite{arnold02uad,riviere08dgm}; dG discretization for periodic boundary conditions is given in \cite{vemaganti07dgm}. The SIPG semi-discretized system of \eqref{varformallencahn} reads as: for a.e. $t \in (0,T]$, find $u_h(\mu)$ in the SIPG finite element space $V_h$ such that
\begin{linenomath*}
\begin{align}\label{discreteallencahn}
(\partial_{t}u_h(\mu), \upsilon_h)_{\Omega} + a_h(\mu;u_h(\mu), \upsilon_h) + (f(u_h(\mu);\mu ),\upsilon_h)_{\Omega}&=0, & \forall \upsilon_h \in V_h,
\end{align}
\end{linenomath*}
with the SIPG bilinear form
\begin{linenomath*}
\begin{equation}\label{bilin}
\begin{aligned}
a_{h}(\mu;u ,\upsilon) =& \sum_{K \in \mathcal{T}_h}\int_K \epsilon \nabla u \cdot\nabla \upsilon - \sum_{E\in E^{0}_{h}}\int_E\left\{\epsilon\nabla u\right\}[\upsilon]ds\\
 & - \sum_{E\in E^{0}_{h}}\int_E\left\{\epsilon\nabla \upsilon \right\}[u]+\sum_{E\in E^{0}_{h}} \frac{\sigma \epsilon}{h_{E}}\int_E [u][\upsilon]ds,
\end{aligned}
\end{equation}
\end{linenomath*}
on a triangulation $\mathcal{T}_h$ with triangular elements $K$ and interior edges $E$ having measure $h_E$. In \eqref{bilin}, $\sigma$ denotes the penalty parameter which should be sufficiently large to ensure the stability of the SIPG scheme \cite{arnold02uad,riviere08dgm}. For easy notation, we omit the explicit dependence of the discrete solution, bilinear form and the non-linear term on the parameter $\mu$. The solution of \eqref{discreteallencahn} is given by
\begin{linenomath*}
\begin{equation*}
u_h(x,t)=\sum^{n_K}_{i=1}\sum^{n_q}_{j=1}u^{i}_{j}(t) \varphi^{i}_{j}(x),
\end{equation*}
\end{linenomath*}
where $\varphi^{i}_{j}(x)$ and $u^{i}_{j}(t)$, $i=1, \ldots, n_k$, $j=1, \ldots, n_q$, are the basis functions of $V_h$ and the unknown coefficients, respectively. The number $n_q$ denotes the local dimension, depending on the order $q$ of the basis functions, on each triangular element and $n_K$ is the number of triangular elements. The unknown coefficients  and basis functions are defined as vectors:
\begin{linenomath*}
\begin{align*}
\bm{u} &:= \bm{u}(t)=(u_{1}^{1}(t), u_{2}^{1}(t), \ldots , u_{n_q}^{n_K}(t))^T :=(u_{1}(t), u_{2}(t), \ldots , u_{\mathcal N}(t))^T,\\
\bm{\varphi} &:= \bm{\varphi}(x)=(\varphi_{1}^{1}(x), \varphi_{2}^{1}(x), \ldots , \varphi_{n_q}^{n_K}(x))^T :=(\varphi_{1}(x), \varphi_{2}(x), \ldots , \varphi_{\mathcal N}(x))^T.
\end{align*}
\end{linenomath*}
Here ${\mathcal N}=n_K\times n_q$ denotes the dG degrees of freedom (DoFs). Then, the SIPG semi-discretized system \eqref{discreteallencahn} leads to the full order model (FOM), in form of a semi-linear system of ordinary differential equations (ODEs):
\begin{linenomath*}
\begin{equation}\label{fom}
\bm{M}\bm{u}_{t} + \bm{A}\bm{u} + \bm{f}(\bm{u})=\bm{0},
\end{equation}
\end{linenomath*}
for the unknown coefficient vector $\bm{u}(t)$,  where $\bm{M}\in\mathbb{R}^{{\mathcal N}\times {\mathcal N}}$ is the mass matrix, $\bm{A}\in\mathbb{R}^{{\mathcal N}\times {\mathcal N}}$ is the stiffness matrix, and $\bm{f}(\bm{u})\in\mathbb{R}^{{\mathcal N}}$ is the non-linear vector of unknown coefficients $\bm{u}$, whose $i$-th entry is given by $\bm{f}_i(\bm{u})=(f(u_{h}(t)),\varphi_i(x))_{\Omega}$, $i=1,\ldots ,\mathcal{N}$.

For the temporal discretization, we consider the uniform partition $0=t_0<t_1<\ldots < t_J=T$ of the time interval $[0,T]$ with the uniform time step-size $\Delta t=t_{n+1}-t_{n}$, $n=0,1,\ldots , J-1$. As the time integrator, we use the AVF method \cite{celledoni12ped,hairer10epv} which preserves the energy decreasing property without restriction of the step size $\Delta t$. The AVF method for a general gradient system $\dot{y} = - \nabla G(y)$ is given  as:
\begin{linenomath*}
\begin{equation*}
y_{n+1} = y_{n} - \Delta t \int _{0}^{1} \nabla G (\tau y_{n+1}+( 1-\tau )y_n)d \tau.
\end{equation*}
\end{linenomath*}

The application of the AVF time integrator to  \eqref{fom} leads to the fully discrete system
\begin{linenomath*}
\begin{equation}\label{fom1}
\bm{M}\bm{u}^{n+1}- \bm{M}\bm{u}^{n} + \frac{\Delta t}{2}\bm{A}(\bm{u}^{n+1} + \bm{u}^{n}) + \Delta t\int^{1}_{0} \bm{f}( \tau \bm{u}^{n+1} + (1-\tau)\bm{u}^{n} ) d\tau =\bm{0}.
\end{equation}
\end{linenomath*}

\subsection{Energy stability of the full order model}
\label{subsec:energyfull}

Now, we prove that the  SIPG-AVF full discretized gradient system is unconditionally energy stable.  The SIPG discretized energy function of the continuous energy $\mathcal{E}(u)$ at a time $t_n = n \Delta t$ is given by:
\begin{linenomath*}
\begin{equation}\label{discreteenergy}
\begin{aligned}
\mathcal{E}_h(u_h^{n}) =& \frac{\epsilon}{2}\left\|\nabla u_h^{n}\right\|^2_{L^2(\Omega)} + (F(u_h^n),1)_{\Omega} \\
& + \sum_{E\in E^{0}_{h}} \left( -(\{\epsilon\nabla u_h^{n}\}, [u_h^{n}])_E + \frac{\sigma\epsilon}{2h_E}([u_h^{n}],[u_h^{n}])_E \right),
\end{aligned}
\end{equation}
\end{linenomath*}
where $u_h^n:=u_h(t_n)\in V_h$. Applying the AVF time integrator to the semi-discrete system \eqref{discreteallencahn} and using the bilinearity of $a_h$, we get for $n=0,1,\ldots , J-1$
\begin{linenomath*}
\begin{align*}
&\frac{1}{\Delta t }(u^{n+1}_{h} - u^{n}_{h}, \upsilon_{h})_{\Omega}  + \frac{1}{2}a_h(u^{n+1}_{h}+u^{n}_{h},\upsilon_{h})\\
 &+ \int_0^1(f(\tau u^{n+1}_{h} + (1-\tau)u^{n}_{h}), \upsilon_{h})_{\Omega} d\tau = 0.
\end{align*}
\end{linenomath*}
Choosing $\upsilon_{h}=u^{n+1}_{h}-u^{n}_{h}$ and using the algebraic identity $(a+b)(a-b)=a^2 - b^2$
\begin{linenomath*}
\begin{equation}\label{eq2avff}
\begin{aligned}
&\frac{1}{ \Delta t}(u^{n+1}_{h} - u^{n}_{h}, u^{n+1}_{h}-u^{n}_{h})_{\Omega} + \frac{1}{2}a_h(u^{n+1}_{h},u^{n+1}_{h}) - \frac{1}{2}a_h(u^{n}_{h},u^{n}_{h}) \\
& + \int_{\Omega}\left[\int_0^1(f(\tau u^{n+1}_{h} + (1-\tau)u^{n}_{h})( u^{n+1}_{h}-u^{n}_{h}) d\tau\right]dx =0.
\end{aligned}
\end{equation}
\end{linenomath*}
By the change of variable $z_h=\tau u^{n+1}_{h} + (1-\tau)u^{n}_{h}$, we get
\begin{linenomath*}
\begin{equation}\label{ftc}
\int_0^1(f(\tau u^{n+1}_{h} + (1-\tau)u^{n}_{h})( u^{n+1}_{h}-u^{n}_{h}) d\tau = \int_{u^{n}_{h}}^{u^{n+1}_{h}}f(z_h) dz= F(u^{n+1}_{h}) - F(u^{n}_{h}).
\end{equation}
\end{linenomath*}
Finally, substituting \eqref{ftc} into \eqref{eq2avff}, using \eqref{bilin} and \eqref{discreteenergy}, we obtain
\begin{linenomath*}
\begin{align*}
\mathcal{E}_h(u^{n+1}_{h}) - \mathcal{E}_h(u^{n}_{h}) &= -\frac{1}{\Delta t}\|u^{n+1}_{h} - u^{n}_{h}\|_{L^2(\Omega)}^2 \leq 0,
\end{align*}
\end{linenomath*}
which implies that $\mathcal{E}_h(u^{n+1}_{h}) \leq  \mathcal{E}_h(u^{n}_{h}) $ for any time step size $\Delta t > 0.$

\section{Model order reduction for gradient systems}
\label{sec:rom}

In this section, we describe the construction of the reduced order model (ROM) and DEIM of the non-linear term for the SIPG discretization. We also show that the reduced solutions using DEIM provides conditional energy stability of the discrete energy function.

\subsection{Reduced order model}
\label{subsec:rom}

The ROM solution $u_{h,r}(x,t)$ of dimension $N\ll {\mathcal  N}$ is formed by approximating the solution $u_h(x,t)$ in a subspace $V_{h,r}\subset V_h$ spanned by a set of $L^2$-orthogonal  basis functions $\{\psi_{i}\}_{i=1}^N$ of dimension $N$, and then projecting onto $V_{h,r}$:
\begin{linenomath*}
\begin{equation}\label{fhn_combrom}
u_h(x,t)\approx u_{h,r}(x,t) = \sum_{i=1}^{N} u_{i,r}(t) \psi_{i}(x) , \qquad (\psi_{i}(x),\psi_{j}(x))_{\Omega}=\delta_{ij},
\end{equation}
\end{linenomath*}
where $\bm{u}_r(t):=(u_{1,r}(t),\ldots ,u_{N,r}(t))^T$ is the coefficient vector of the reduced solution.  Then, the SIPG weak formulation for ROM reads as: for a.e. $t\in (0,T]$, find $u_{h,r}(x,t)\in V_{h,r}$ such that
\begin{linenomath*}
\begin{align}\label{rom_weak}
(\partial_{t}u_{h,r}, \upsilon_{h,r})_{\Omega} + a_{h}(u_{h,r},\upsilon_{h,r}) + (f(u_{h,r}),\upsilon_{h,r})_{\Omega}&=0, & \forall \upsilon_{h,r}\in V_{h,r}
\end{align}
\end{linenomath*}

Since the reduced basis functions $\{\psi_{i}\}_{i=1}^N\subset V_{h,r}$ also belong to the space $V_h$, they can be expanded by finite element basis functions  $\{\varphi_i(x)\}_{i=1}^{\mathcal{N}}$ as:
\begin{equation}\label{fhn_combpod}
\psi_{i}(x) = \sum_{j=1}^{\mathcal{N}} \Psi_{j,i} \varphi_j(x), \qquad \Psi_{\cdot ,i}^T\bm{M}\Psi_{\cdot ,j}=\delta_{ij}.
\end{equation}
The coefficient vectors of the  reduced basis function are collected in the columns of the matrix $\bm{\Psi} =[\Psi_{\cdot ,1}, \ldots ,\Psi_{\cdot ,N} ]\in\mathbb{R}^{{\mathcal{N}}\times N}$. The  coefficient vectors of FOM and ROM solutions are related by $\bm{u} \approx \bm{\Psi}\bm{u}_r$. Substituting this relation together with \eqref{fhn_combrom} and \eqref{fhn_combpod} into the system \eqref{rom_weak}, we obtain for the unknown coefficient vectors the reduced semi-discrete ODE system:
\begin{equation}\label{rom}
\partial_t\bm{u}_r  + \bm{A}_r\bm{u}_r + \bm{f}_r(\bm{u}_r) = \bm{0},
\end{equation}
with the reduced stiffness matrix $\bm{A}_r=\bm{\Psi}^T\bm{A}\bm{\Psi}$ and the reduced non-linear vector $\bm{f}_r(\bm{u}_r)=\bm{\Psi}^T\bm{f}(\bm{\Psi}\bm{u}_r)$. The construction of the reduced basis functions $\{\psi_{i}\}_{i=1}^N$ is discussed in Section~\ref{subsec:sampling}.

\subsection{Discrete empirical interpolation method (DEIM)}
\label{subsec:deim}

Although the dimension of the reduced system \eqref{rom} is small, $N\ll {\mathcal N }$, the computation of the reduced non-linear vector $\bm{f}_r(\bm{u}_r)=\bm{\Psi}^T\bm{f}(\bm{\Psi}\bm{u}_r)$ still depends on the dimension ${\mathcal N }$ of the full system. In order to reduce the online computational cost, we apply the DEIM \cite{chaturantabut10nmr} to approximate the non-linear vector $\bm{f}(\bm{\Psi}\bm{u}_r)\in\mathbb{R}^{\mathcal{N}}$ from a $M\ll {\mathcal N}$ dimensional subspace spanned by non-linear vectors $\bm{f}(\bm{\Psi}\bm{u}_r(t_n))$, $n=1,\ldots , J$. Let $M\ll {\mathcal N}$ orthonormal  basis functions $\{W_i\}_{i=1}^M$ are given. We set the matrix $\bm{W}:=[W_1,\ldots , W_M]\in\mathbb{R}^{\mathcal{N}\times M}$ (the functions $W_i$ are computed successively during the greedy iteration in EIM, whereas here, in DEIM, the functions $W_i$ are computed priori by POD and then they are used in the greedy iteration). Then, we can use the approximation $\bm{f}(\bm{\Psi}\bm{u}_r) \approx \bm{Q}\bm{f}_{m}(\bm{\Psi}\bm{u}_r)$, where $\bm{f}_m(\bm{\Psi}\bm{u}_r)=\bm{P}^T \bm{f}(\bm{\Psi}\bm{u}_r)\in\mathbb{R}^{M}$ and the matrix $\bm{Q}=\bm{W}(\bm{P}^T\bm{W})^{-1}\in\mathbb{R}^{\mathcal{N}\times M}$ is precomputable. For the details of the computation of the reduced non-linear vectors we refer to the greedy DEIM algorithm  \cite{chaturantabut10nmr}.  For continuous finite element and finite volume discretizations, the number of flops for the computation of bilinear form and nonlinear term depends on the maximum number of neighbor cells \cite{Drohmann12}. In the case of dG discretization, due to its local nature, it depends only on the number of nodes in the local cells. For instance, in the case of SIPG with linear elements ($n_q=3$), for each degree of freedom, integrals have to be computed  on a single triangular element \cite{karasozen15spi}, whereas in the case of continuous finite elements, integral computations on $6$ neighbor cells are needed \cite{Heinkenschloss14}. Since the AVF method is an implicit time integrator, at each time step, a non-linear system of equations has to be solved by Newton's method. The reduced Jacobian has a diagonal block structure for the SIPG discretization, which is easily invertible \cite{karasozen15spi}, and requires $O(n_qM)$ operations with DEIM.

\subsection{Energy stability of the reduced solution}
\label{subsec:energyrom}

The energy stability of the DEIM reduced order model is proved in the same way as for the FOM. Applying the AVF time integrator to the semi-discrete system \eqref{rom_weak}, choosing $\upsilon_{h,r}=u^{n+1}_{h,r}-u^{n}_{h,r}$, and using the algebraic identity $(a+b)(a-b)=a^2 - b^2$ and the bilinearity of $a_h$, we obtain
\begin{linenomath*}
\begin{equation}\label{eq2avffr}
\begin{aligned}
\frac{1}{ \Delta t}(u^{n+1}_{h,r} - u^{n}_{h,r}, u^{n+1}_{h,r}-u^{n}_{h,r})_{\Omega} + \frac{1}{2}a_h(u^{n+1}_{h,r},u^{n+1}_{h,r}) - \frac{1}{2}a_h(u^{n}_{h,r},u^{n}_{h,r}) \\
+ \int_0^1\left[\int_{\Omega} f(\tau u^{n+1}_{h,r} + (1-\tau)u^{n}_{h,r})( u^{n+1}_{h,r}-u^{n}_{h,r}) dx\right]d\tau =0.
\end{aligned}
\end{equation}
\end{linenomath*}
Let  us set the averaged reduced solution $z_{h,r}=\tau u^{n+1}_{h,r} + (1-\tau)u^{n}_{h,r}$, and the averaged coefficient vector $\bm{z}_r=\tau \bm{u}_r^{n+1} + (1-\tau)\bm{u}_r^{n}$ of the reduced system. Then, from the integral term in \eqref{eq2avffr}:
\begin{linenomath*}
\begin{align*}
\int_0^1\left[\int_{\Omega} f(z_{h,r})( u^{n+1}_{h,r}-u^{n}_{h,r}) dx\right]d\tau = & \int_0^1\left[\sum_{i=1}^N(u^{n+1}_{i,r}(t)-u^{n}_{i,r}(t))\overbrace{\int_{\Omega} f(z_{h,r})\psi_i(x)dx}^{(\bm{\Psi}^T\bm{f}(\bm{\Psi}\bm{z}_r))_i}\right]d\tau \\
= & \int_0^1 (\bm{\Psi}(\bm{u}_r^{n+1}-\bm{u}_r^{n}))^T \bm{f}(\bm{\Psi}\bm{z}_r))d\tau.
\end{align*}
\end{linenomath*}
On the left hand side of \eqref{eq2avffr}, using the DEIM approximation $\bm{f}(\bm{\Psi}\bm{z}_r)\approx \bm{Q}\bm{f}_{m}(\bm{\Psi}\bm{z}_r)$, adding and subtracting the term $\int_0^1\left[\int_{\Omega} f(z_{h,r})( u^{n+1}_{h,r}-u^{n}_{h,r}) dx\right]d\tau$, and using the integral mean theorem, we obtain
\begin{linenomath*}
\begin{align*}
\mathcal{E}_h(u^{n+1}_{h,r}) - \mathcal{E}_h(u^{n}_{h,r}) = &-\frac{1}{\Delta t}\|u^{n+1}_{h,r} - u^{n}_{h,r}\|_{L^2(\Omega)}^2\\
&+ \int_0^1 (\bm{\Psi}(\bm{u}_r^{n+1}-\bm{u}_r^{n}))^T( \bm{f}(\bm{\Psi}\bm{z}_r) - \bm{Q}\bm{f}_{m}(\bm{\Psi}\bm{z}_r) )d\tau\\
=& -\frac{1}{\Delta t}\|u^{n+1}_{h,r} - u^{n}_{h,r}\|_{L^2(\Omega)}^2 \\
&+ \langle \bm{\Psi}(\bm{u}_r^{n+1}-\bm{u}_r^{n}) ,  \bm{f}(\bm{\Psi}\bm{\tilde{z}}_r^n) - \bm{Q}\bm{f}_{m}(\bm{\Psi}\bm{\tilde{z}}_r^n) \rangle,
\end{align*}
\end{linenomath*}
for some $\bm{\tilde{z}}_r^n$ between $\bm{u}_r^{n}$ and $\bm{u}_r^{n+1}$, and $\langle \cdot , \cdot \rangle$ denoting the Euclidean inner product. Applying the Cauchy-Schwarz inequality, we get
\begin{linenomath*}
\begin{equation*}
\langle \bm{\Psi}(\bm{u}_r^{n+1}-\bm{u}_r^{n}) ,  \bm{f}(\bm{\Psi}\bm{\tilde{z}}_r^n) - \bm{Q}\bm{f}_{m}(\bm{\Psi}\bm{\tilde{z}}_r^n) \rangle \leq \|\bm{\Psi}(\bm{u}_r^{n+1}-\bm{u}_r^{n})\|_2\|\bm{f}(\bm{\Psi}\bm{\tilde{z}}_r^n) - \bm{Q}\bm{f}_{m}(\bm{\Psi}\bm{\tilde{z}}_r^n)\|_2.
\end{equation*}
\end{linenomath*}
Using the a priori error bound \cite{chaturantabut10nmr,Chaturantabut12a,Heinkenschloss14}, we have
\begin{linenomath*}
\begin{equation*}
\|\bm{f}(\bm{\Psi}\bm{\tilde{z}}_r^n) - \bm{Q}\bm{f}_{m}(\bm{\Psi}\bm{\tilde{z}}_r^n)\|_2  \leq \| (\bm{P}^T\bm{W} )^{-1}||_2 \|(\bm{I} -\bm{W}\bm{W}^T)\bm{f}(\bm{\Psi}\bm{\tilde{z}}_r^n)||_2.
\end{equation*}
\end{linenomath*}

Using the equivalent weighted-Euclidean inner product form of the $L^2$-norm on the reduced space $V_{h,r}$, we have
\begin{linenomath*}
\begin{align*}
\|u^{n+1}_{h,r} - u^{n}_{h,r}\|_{L^2(\Omega)}^2 =& (\bm{u}_r^{n+1}-\bm{u}_r^{n})^T\bm{M}_r(\bm{u}_r^{n+1}-\bm{u}_r^{n}) = (\bm{u}_r^{n+1}-\bm{u}_r^{n})^T\bm{\Psi}^T\bm{M}\bm{\Psi}(\bm{u}_r^{n+1}-\bm{u}_r^{n})\\
=& (\bm{u}_r^{n+1}-\bm{u}_r^{n})^T\bm{\Psi}^T\bm{R}^T\bm{R}\bm{\Psi}(\bm{u}_r^{n+1}-\bm{u}_r^{n}) = \| \bm{R}\bm{\Psi}(\bm{u}_r^{n+1}-\bm{u}_r^{n}) \|_2^2,
\end{align*}
\end{linenomath*}
where $\bm{M}_r=\bm{\Psi}^T\bm{M}\bm{\Psi}$ is the reduced mass matrix (indeed it is the identity matrix, $\bm{\Psi}$ is $\bm{M}$-orthogonal), and $\bm{R}$ is the Cholesky factor of the mass matrix $\bm{M}$ (i.e. $\bm{M}=\bm{R}^T\bm{R}$). Thus, we get the identity
\begin{linenomath*}
\begin{align*}
\| \bm{\Psi}(\bm{u}_r^{n+1}-\bm{u}_r^{n}) \|_2 =& \| \bm{R}^{-1}\bm{R}\bm{\Psi}(\bm{u}_r^{n+1}-\bm{u}_r^{n}) \|_2 \leq  \| \bm{R}^{-1}\|_2 \|\bm{R}\bm{\Psi}(\bm{u}_r^{n+1}-\bm{u}_r^{n}) \|_2\\
=& \| \bm{R}^{-1}\|_2 \|u^{n+1}_{h,r} - u^{n}_{h,r}\|_{L^2(\Omega)}.
\end{align*}
\end{linenomath*}
Using the above identity, we obtain for the energy difference:
\begin{linenomath*}
\begin{equation}\label{energydiff}
\begin{aligned}
\mathcal{E}_h(u^{n+1}_{h,r}) - \mathcal{E}_h(u^{n}_{h,r}) \leq & \|u^{n+1}_{h,r} - u^{n}_{h,r}\|_{L^2(\Omega)}^2 \; \Cross \\
  &  \left(  -\frac{1}{\Delta t} + \frac{\| \bm{R}^{-1}\|_2 \| (\bm{P}^T\bm{W} )^{-1}||_2 \|(\bm{I} -\bm{W}\bm{W}^T)\bm{f}(\bm{\Psi}\bm{\tilde{z}}_r^n)||_2}{\|u^{n+1}_{h,r} - u^{n}_{h,r}\|_{L^2(\Omega)} } \right).
\end{aligned}
\end{equation}
\end{linenomath*}
The ROM satisfies the energy decrease property $\mathcal{E}_h(u^{n+1}_{h,r}) \leq \mathcal{E}_h(u^{n}_{h,r})$  when the right hand side of \eqref{energydiff} is non-positive, i.e., if the time-step size is bounded as
\begin{linenomath*}
\begin{align}\label{upperbound}
\Delta t\leq \frac{\|u^{n+1}_{h,r} - u^{n}_{h,r}\|_{L^2(\Omega)}}{\| \bm{R}^{-1}\|_2 \| (\bm{P}^T\bm{W} )^{-1}||_2 \|(\bm{I} -\bm{W}\bm{W}^T)\bm{f}(\bm{\Psi}\bm{\tilde{z}}_r^n)||_2} .
\end{align}
\end{linenomath*}

The columns of the matrix $\bm{W}$ are orthonormal and $\|(\bm{P}^T\bm{W} )^{-1}\|_2$  and $\|\bm{R}^{-1}\|_2$ are of moderate size. They differ in the numerical tests between $10-30$, and $30-60$, respectively. The upper bound for $\Delta t$ in \eqref{upperbound} for each time step can be extended for all time steps to the following global upper bound:
\begin{linenomath*}
\begin{align}\label{globalupperbound}
\Delta t \leq \frac{\|u_{h,r}\|_{m,L^2(\Omega)}}{\| \bm{R}^{-1}\|_2 \| (\bm{P}^T\bm{W} )^{-1}||_2 \|\bm{f}(\bm{\Psi}\bm{\tilde{z}}_r)||_{M,2}},
\end{align}
\end{linenomath*}
where
\begin{linenomath*}
\begin{align*}
\|u_{h,r}\|_{m,L^2(\Omega)} &= \min_{1\leq n\leq J-1} \|u^{n+1}_{h,r} - u^{n}_{h,r}\|_{L^2(\Omega)},\\
\|\bm{f}(\bm{\Psi}\bm{\tilde{z}}_r)||_{M,2}  &= \max_{1\leq n\leq J-1} \|(\bm{I} -\bm{W}\bm{W}^T)\bm{f}(\bm{\Psi}\bm{\tilde{z}}_r^n)||_2.
\end{align*}
\end{linenomath*}

The global upper bound for the time step-size $\Delta t$ in the right hand side of \eqref{globalupperbound}, is a sufficiently large number so that we can choose $\Delta t$ sufficiently large in the numerical examples in Section~\ref{numerical}. Hence, the DEIM reduced energy decreases almost unconditionally for large time step-size.

\subsection{POD Greedy Adaptive Sampling}
\label{subsec:sampling}

For an efficient offline-online computation of the reduced basis functions $\{\psi_i\}_{i=1}^N$, several greedy sampling algorithms are developed for finite difference, finite element and finite volume methods (see for example \cite{Drohmann12,Grepl13,Wirtz14}). In this section we will describe the POD-greedy sampling procedure for the SIPG discretized Allen-Cahn equation.  Let $\bm{u}_{\mu^*}$ denotes the solution vector of the FOM related to a parameter value $\mu^*$. Let us also denote by $POD_{\bm{X}}(\bm{B},k)$ the operator on the Euclidean space with $\bm{X}$-weighted inner product ($\bm{X}$ is a positive definite matrix), and giving $k$ POD basis functions of the matrix $\bm{B}$, related to the first $k$ largest singular values. For instance, for a parameter value $\mu^*$, the $\bm{M}$-orthonormal reduced modes $\{\Psi_i\}_{i=1}^N$, coefficient vectors of the $L^2$-orthonormal reduced basis functions $\{\psi_i\}_{i=1}^N$, may be computed through the generalized singular value decomposition \cite{ karasozen15spi} of the snapshot matrix $[\bm{u}_{\mu^*}^{n_1},\ldots , \bm{u}_{\mu^*}^J]\in\mathbb{R}^{{\mathcal N}\times J}$ by $\{\Psi_1,\ldots , \Psi_N\}=POD_{\bm{M}}([\bm{u}_{\mu^*}^{n_1},\ldots , \bm{u}_{\mu^*}^J],N)$.

Because the computation of the POD modes for time dependent parametrized PDEs is computationally demanding, we develop  an adaptive  POD-greedy (PODG) algorithm, where we perform a greedy search among a parameter space  $\mathcal{M}$ \cite{Grepl13}. The algorithm  starts by selecting an initial parameter set $\mathcal{M}_0=\{\mu^*\}$, where $\mu^*$ belongs to a training set $\mathcal{M}_{train}=\{\mu_1,\ldots ,\mu_{n_s}\}\subset\mathcal{M}$, and an empty reduced space $V_{h,r}^{0}=\{0\}$. At the $k$-th greedy iteration, we determine the parameter $\mu^*\in \mathcal{M}_{train}$ for which an error indicator $\Delta_k(\mu^*)$ is related to the reduced system \eqref{rom} on $V_{h,r}^{k-1}$. Then, we extend the reduced space $V_{h,r}^{k-1}$ by adding a single POD mode corresponding to the dominant singular value of  $\bm{e}_{\mu^*}:=[\bm{e}_{\mu^*}^1,\ldots ,\bm{e}_{\mu^*}^J]$, where $\bm{e}_{\mu^*}^n=\bm{u}_{\mu^*}^n-\text{Proj}_{V_{h,r}^{k-1}}\bm{u}_{\mu^*}^n$ is the projection error on $V_{h,r}^{k-1}$. Here, the  single POD mode is computed by the operator $POD_{\bm{M}}(\bm{e}_{\mu^*},1)$. We stop the greedy iteration either until a predefined maximum number $N_{max}$ is reached, or the  error indicator $\Delta_k(\mu^*)$ is below a prescribed tolerance $TOL_G$. We use as an error  indicator the residual-based a-posteriori indicator
\begin{linenomath*}
\begin{equation*}
\Delta_k(\mu) = \left( \Delta t\sum_{n=1}^J \|R_{h}(u_{r,\mu}^n)\|_{H^{-1}} \right)^{1/2},
\end{equation*}
\end{linenomath*}
where $\|\cdot\|_{H^{-1}}$ is the dual norm on $H^1$, and $R_h(u_{r,\mu}^n)$ denotes the residual of the reduced system \eqref{rom_weak} on the $n$-th time level after time discretization.

In addition, at each greedy iteration, the computation of the error indicators $\Delta_k(\mu_i)$, $i=1,\ldots , n_s$, requires the solution of the reduced system \eqref{rom} for several times. For an efficient offline/online decomposition, we make use of the affine dependence of the bilinear form $a_h$ on the parameter $\epsilon$. Thus the stiffness matrix $\bm{A}^1$ related to the bilinear form $a_h(1;u,v)$ is computed in the offline stage only once. In the online stage, $\bm{A}_r=\epsilon\bm{\Psi}^T\bm{A}^1\bm{\Psi}$ is computed without an additional cost. On the other hand, the non-linear term related to the logarithmic potential does not depend affinely on the parameter $\theta$, therefore in the online stage, the non-linear vector in \eqref{rom} is approximated using DEIM, which requires a quite small number of operations, $M\ll {\mathcal N }$ for nonlinear vector and $n_qM$ for Jacobian computation, and the matrix $\bm{Q}$ in the DEIM approximation is computed only once in the offline stage. Moreover, in the POD-Greedy basis computation, we also use the DEIM approximation ("Inner DEIM" in the  Algorithm~\ref{greedy_alg}). Related to a certain parameter value $\mu^*$, the (temporary) inner DEIM basis functions $\{W_1,\ldots , W_M\}$ are computed through the POD of the snapshot ensemble $\mathcal{F}_*:=[\bm{f}^1_{\mu^*}, \ldots ,\bm{f}^J_{\mu^*}]$  of the non-linear vectors; i.e. $\{W_1,\ldots , W_M\}:=POD_{\bm{I}}(\mathcal{F}_*,M)$, where $\bm{I}$ is the identity matrix. In order to minimize the error induced by DEIM, we take a sufficiently large $M=M_{max,*}\leq rank(\mathcal{F}_*)$ at each greedy iteration. Finally, we construct the (final)  outer DEIM basis functions by applying the POD to the snapshot matrix collecting all the snapshots for the parameter values $\mu\in \mathcal{M}_{N}$, which are stored in the POD-greedy algorithm.

\begin{algorithm}
\caption{POD-greedy algorithm}
\textbf{Input:} Samples $\mathcal{M}_{train}=\{\mu_i\}$, $|\mathcal{M}_{train}|=n_s$, tolerance $TOL_G$
 \\
\textbf{Output:} $V_{h,r}^{N}:=\text{span}\{\Psi_1,\ldots ,\Psi_N\}$, $\bm{W}:=\text{span}\{W_1,\ldots ,W_M\}$
 \\
\begin{algorithmic}
\STATE  $\mathcal{M}_0:=\{\mu_1\}$, $\mu^*=\mu_1$, $V_{h,r}^{0}:=\{0\}$, $N=1$\\

\WHILE{$N\leq N_{max}$}
   \STATE  compute $\bm{u}^{n}_{\mu^*}$, $n=1,2,\ldots ,J$\\
	 \STATE  set $\bm{e}_{\mu^*}^n=\bm{u}_{\mu^*}^{n} - \text{Proj}_{V_{h,r}^{N-1}}\bm{u}_{\mu^*}^{n}$, $n=1,2,\ldots ,J$ \\
	 \STATE  $V_{h,r}^{N} \; \longleftarrow \; V_{h,r}^{N-1} \cup \text{POD}_{\bm{M}}(\{\bm{e}_{\mu^*}^1,\ldots , \bm{e}_{\mu^*}^J\},1)$\\
	 \FOR{$i=1$ \TO $n_s$}
			\STATE 	$\{W_1,\ldots ,W_{M_{max,i}}\}=\text{POD}_{\bm{I}}(\{\bm{f}_{\mu_i}^1,\ldots , \bm{f}_{\mu_i}^J\},M_{max,i})$ \quad (Inner DEIM Basis)\\
			\STATE solve reduced system on $V_{h,r}^{N}$ using DEIM approximation\\
			\STATE calculate error indicator $\Delta_N(\mu_i)$
	 \ENDFOR
	\STATE $\mu^*=\underset{\mu\in \{\mu_1,\ldots ,\mu_{n_s}\}}{\tt{argmax}} \Delta_N(\mu)$\\
	\IF{ $\Delta_{N}(\mu^*)\leq TOL_G$ }
		   \STATE $N_{max}=N$\\
			 \STATE break \\
	\ENDIF\\	
	\STATE $\mathcal{M}_N:=\mathcal{M}_{N-1}\cup \{\mu^*\}$\\
	\STATE $N \; \longleftarrow \; N+1$\\
\ENDWHILE
\end{algorithmic}

\begin{algorithmic}
\STATE $\mathcal{F} \; \longleftarrow \; [\mathcal{F}_1,\ldots , \mathcal{F}_N]$, \; $\mathcal{F}_i=[\bm{f}_{\mu_i}^1,\ldots , \bm{f}_{\mu_i}^J]$,\; $\mu_i\in \mathcal{M}_N$,\; $i=1,\ldots ,N$ \\
\STATE $\{W_1,\ldots ,W_M\}=\text{POD}_{\bm{I}}(\mathcal{F},M)$  \quad (Outer DEIM Basis)
\end{algorithmic}
\label{greedy_alg}
\end{algorithm}

The solutions exhibit larger gradients  for the sharp interface limit when $\epsilon \rightarrow 0$. In  this case the FOM solutions are computed on spatially non-uniform grids using moving mesh methods \cite{Shen09} or  adaptive finite elements \cite{Feng05,zhang09nsd}. Using space adaptive methods, more points are located at the sharp interface in order to resolve the steep gradients. Compared to the fixed meshes, space-adaptive meshes require less degrees of freedom, which also would increase the speed-up of the ROMs. On the other hand, the dimension of the snapshots changes at each time step for space-adaptive meshes in contrast to the fixed size of the snapshots for fixed meshes. In  order to deal with this problem, a common discretization space can be formed.  But this space would be relatively high dimensional for sharp interface problems with locally varying features. In \cite{Ullmann16} formation of the fixed common discretization space is avoided without interpolating the snapshots. But, then the error of the POD reduced solutions do not satisfy the Galerkin orthogonality to the reduced space created by the adaptive snapshots. An error analysis of the POD Galerkin method for linear elliptic problems is performed in  \cite{Ullmann16} and applicability of the approach is tested for 2D linear convection problems and time dependent Burger's equation. We also mention that spatially adaptive ROMs are studied for adaptive wavelets in \cite{Urban15} and for adaptive mixed finite elements in \cite{Yano16}.

\section{Numerical results}
\label{numerical}

In this section we give two numerical tests to demonstrate the effectiveness of the ROM. In the PODG algorithm, we set the tolerance $TOL_G=10^{-3}$ and the maximum number of PODG basis functions $N_{max}=20$. For the selected parameter values $\mu^*$ in the greedy algorithm, the FOMs  are solved using linear dG elements with uniform spatial mesh size $h:=\Delta x_1=\Delta x_2$. The average number of Newton iterations was one for solving the nonlinear equations (9) at each time step.  In all examples, we present the $L^2(0,T;L^2(\Omega))$ errors of the difference between the FOM and ROM solutions, and $L^{\infty}(0,T)$ errors of the difference between the discrete energies, i.e. the maximum error among the discrete time instances.

\subsection{Allen-Cahn equation with quartic potential functional} 
\label{acq}

We consider the 2D Allen-Cahn equation in \cite{Guo16} with a quartic potential functional \eqref{quartic}, so the bistable non-linear function $f(u)=u^3-u$, under homogeneous Neumann boundary conditions in the  spatial domain $\Omega =[0,1]^2$ and in the time interval $t\in [0,1]$. The spatial and temporal step sizes are taken as $h = 0.015$ and $\Delta t=0.01$, respectively. The initial condition is
\begin{linenomath*}
\begin{equation*}
u(x,0) = \text{tanh}\left( \frac{ 0.25 - \sqrt{ (x_1-0.5)^2 + (x_2-0.5)^2} }{\sqrt{2}\epsilon} \right).
\end{equation*}
\end{linenomath*}

\begin{figure}[htb!]
\centering
\subfloat{\includegraphics[width=0.45\textwidth]{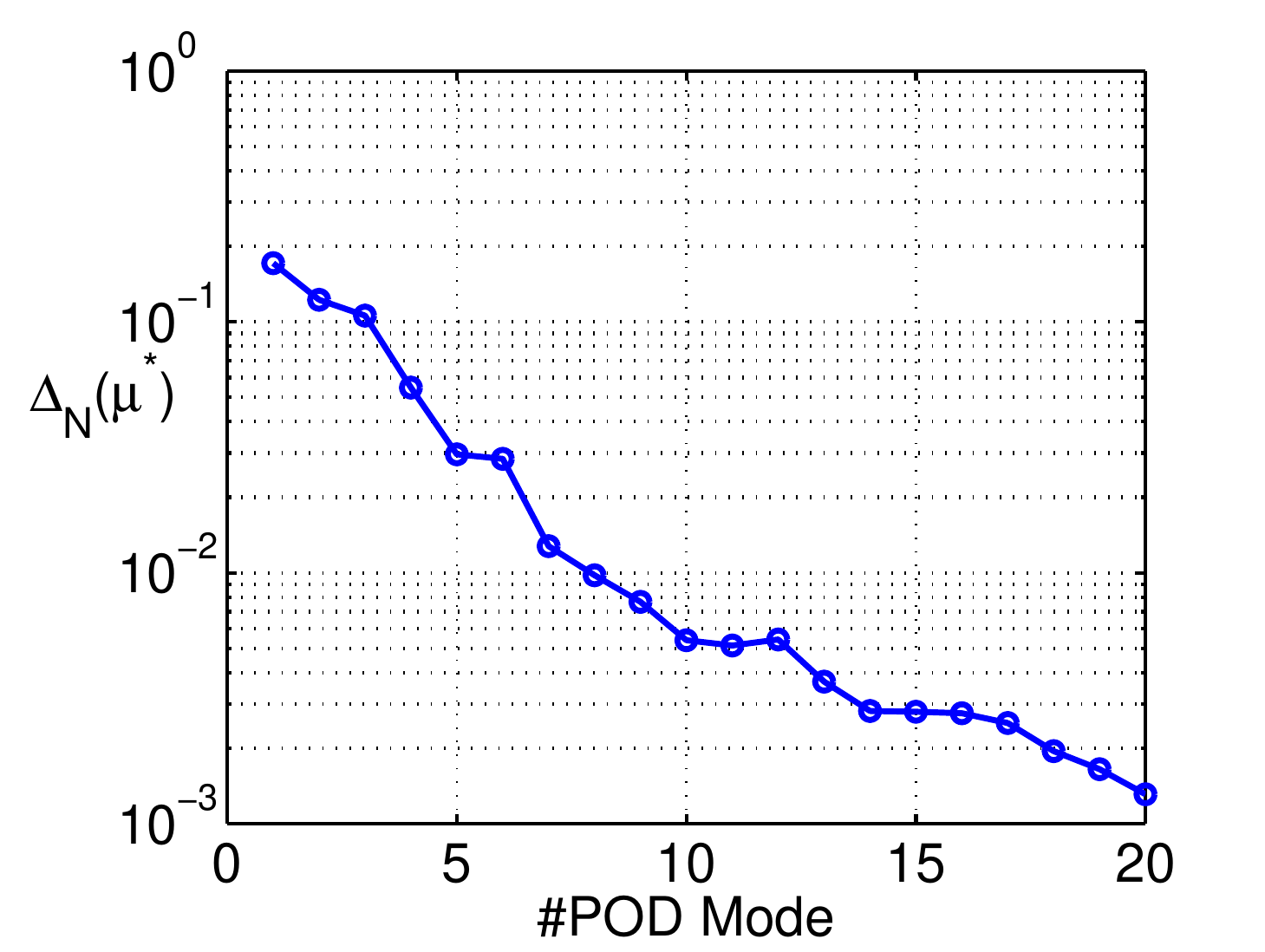}}
\subfloat{\includegraphics[width=0.45\textwidth]{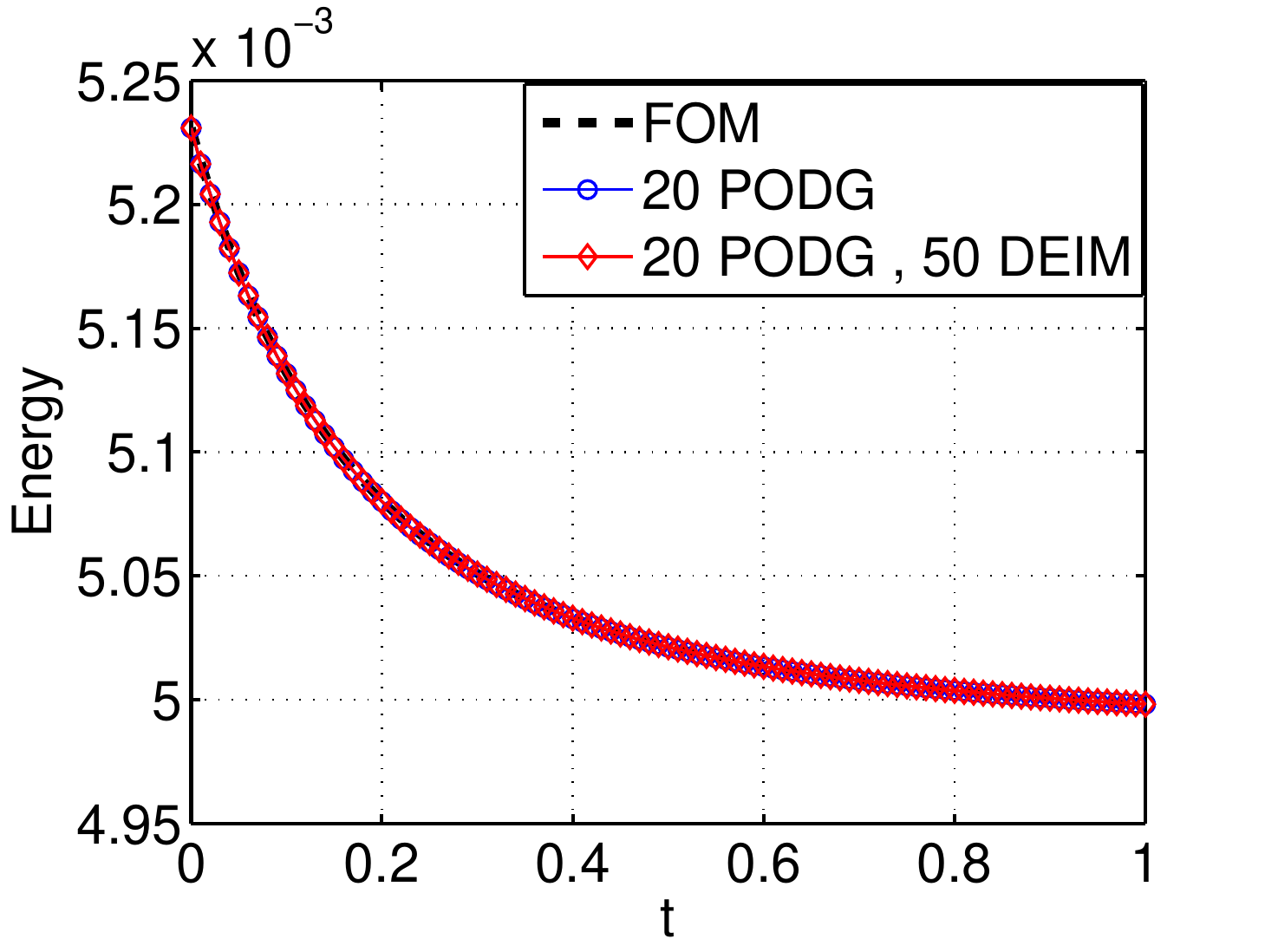}}
\caption{Allen-Cahn  with quartic potential: Error indicator vs POD modes (left) and decrease of energies for $\mu =200$ (right).\label{ac_error}}
\end{figure}

\begin{figure}[htb!]
\centering
\subfloat{\includegraphics[width=0.45\textwidth]{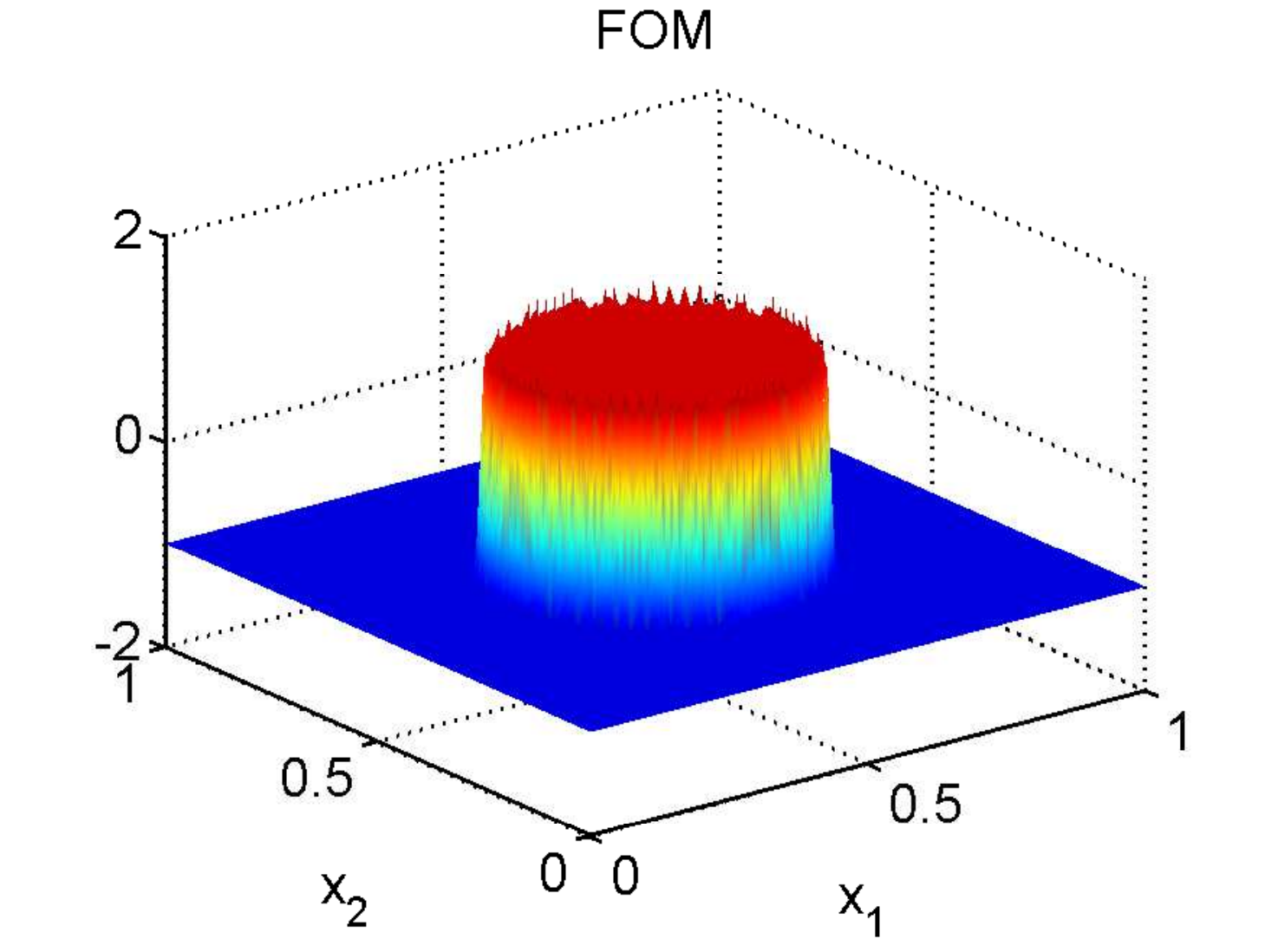}}
\subfloat{\includegraphics[width=0.45\textwidth]{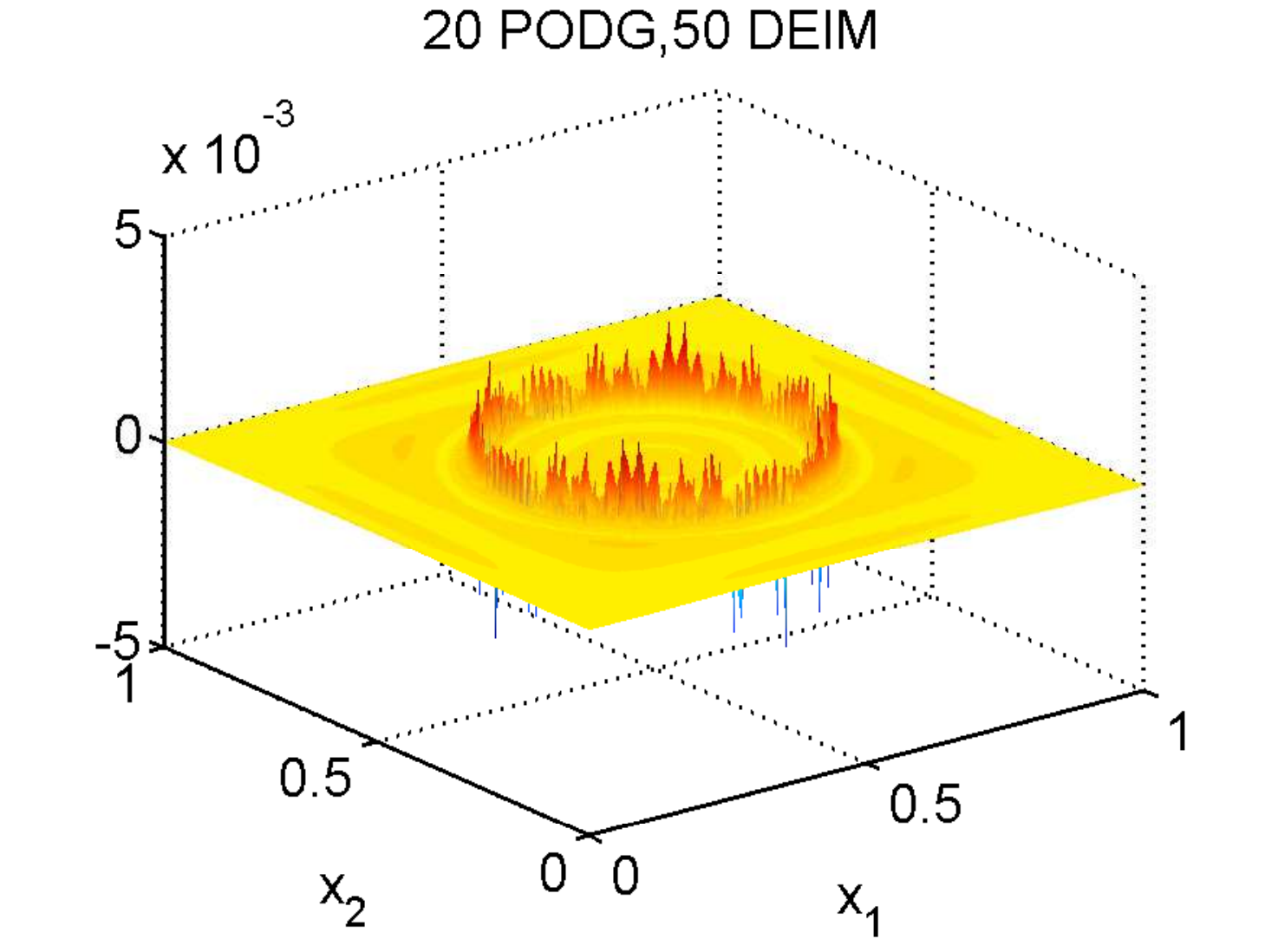}}
\caption{Allen-Cahn  with quartic potential: FOM profile (left) and error plot between FOM and PODG-DEIM (right) solutions, at the final time for $\mu =200$. \label{ac_plot}}
\end{figure}

The FOM becomes stiff for smaller $\epsilon$. The training set for $\mu=1/\epsilon$ is chosen by Clenshaw-Curtis points \cite{Clenshaw60} using more points in direction of larger
$\mu$, or smaller $\epsilon$, respectively:
$${\mathcal M}_{train} =\{10.00,24.78,67.32,132.5,212.46,297.54,377.5,442.68,485.22,500.00\}.$$ The decrease of the error indicator and energy decreases are given in Fig.~\ref{ac_error}. The error plot between FOM and PODG-DEIM solutions at the final time in Fig.~\ref{ac_plot} shows that the dynamic of the system can be captured efficiently by 20 PODG and 50 DEIM modes, which shows that FOM and ROM  solutions of the Allen-Cahn equation is robust with respect to $\epsilon$.

\subsection{Allen-Cahn equation with logarithmic potential functional} 
\label{acl}

Our second example is the 2D Allen-Cahn equation in \cite{Shen16}, obtained here through a scale of the system \eqref{allencahn} by a factor $\sqrt{\epsilon}/2$. We consider the system under periodic boundary conditions, and with a non-convex logarithmic potential \eqref{logarithmic}. We work on $\Omega =[0,2\pi]^2$ with the terminal time $T=1$. For the mesh sizes, we take $\Delta t=0.01$ and $h\approx 0.015$. We accept the initial condition as $u(x,0)= 0.05(2\times \texttt{rand} - 1)$, where the term $\texttt{rand}$ stands for a random number in $[0,1]$.

\begin{figure}[htb!]
\centering
\subfloat{\includegraphics[width=0.45\textwidth]{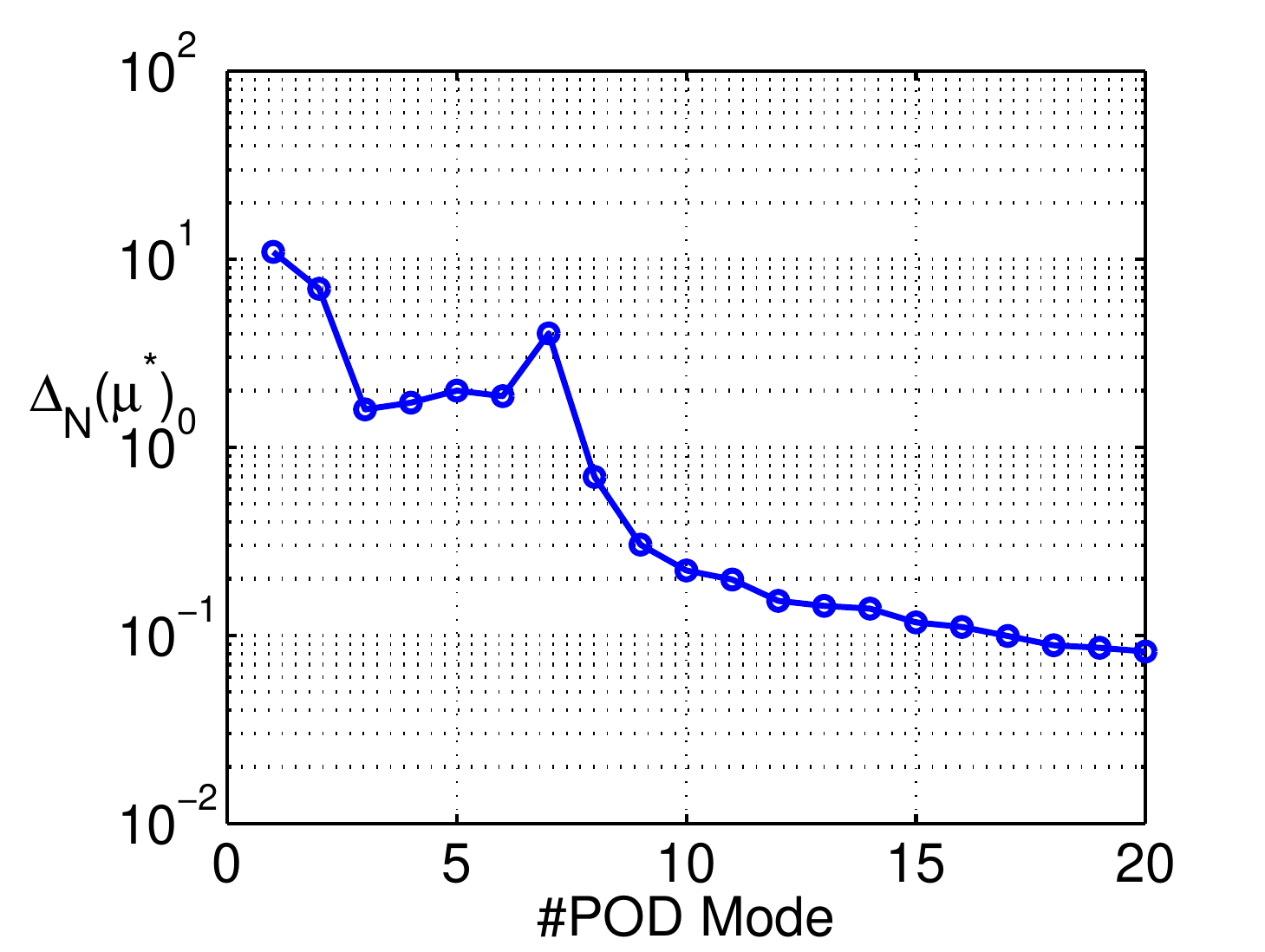}}
\subfloat{\includegraphics[width=0.45\textwidth]{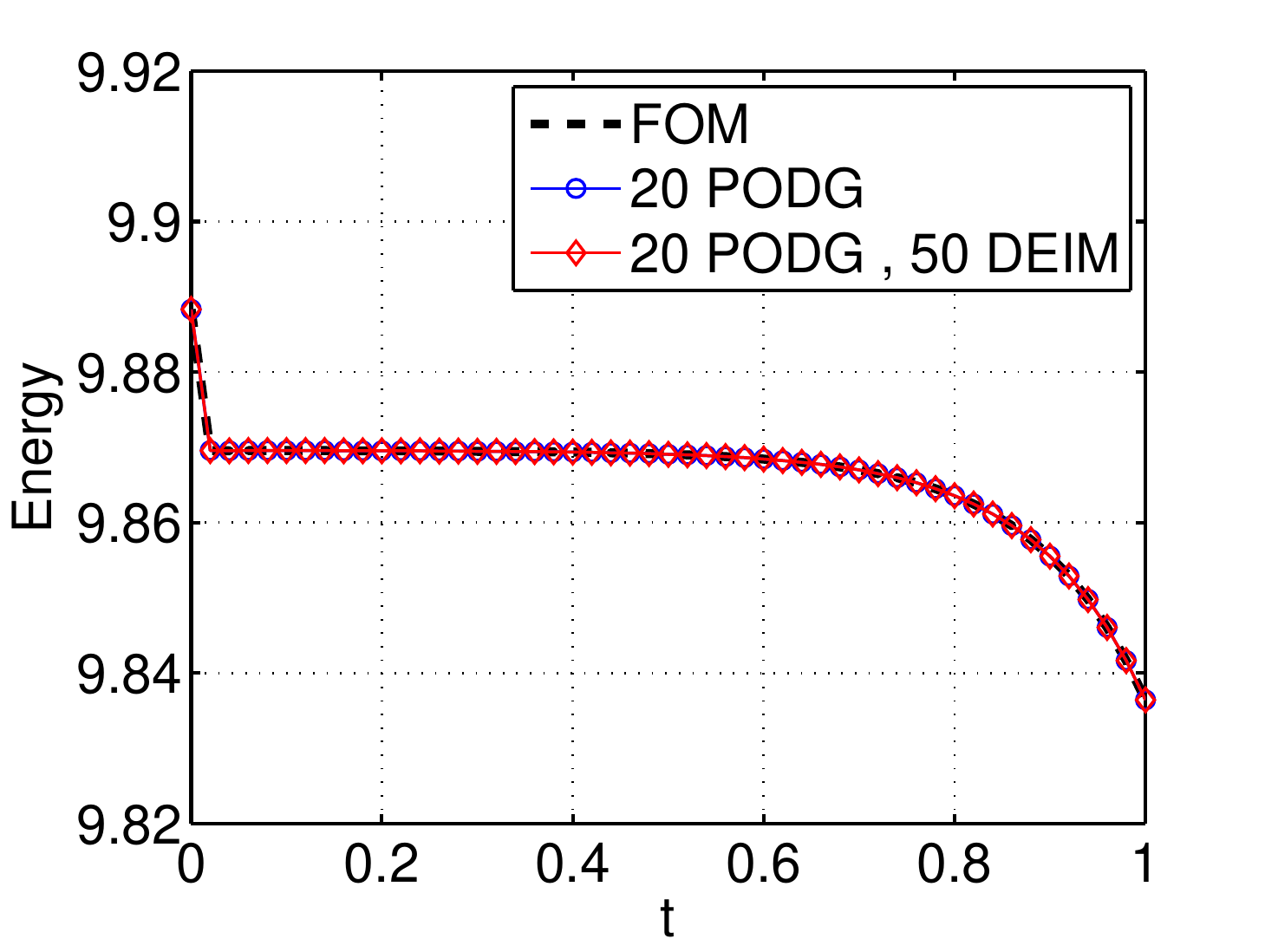}}
\caption{Allen-Cahn with logarithmic potential: Error indicator vs POD modes (left) and decrease of energies for $\theta =0.10$ (right).\label{ac_log_error}}
\end{figure}

\begin{figure}[htb!]
\centering
\subfloat{\includegraphics[width=0.45\textwidth]{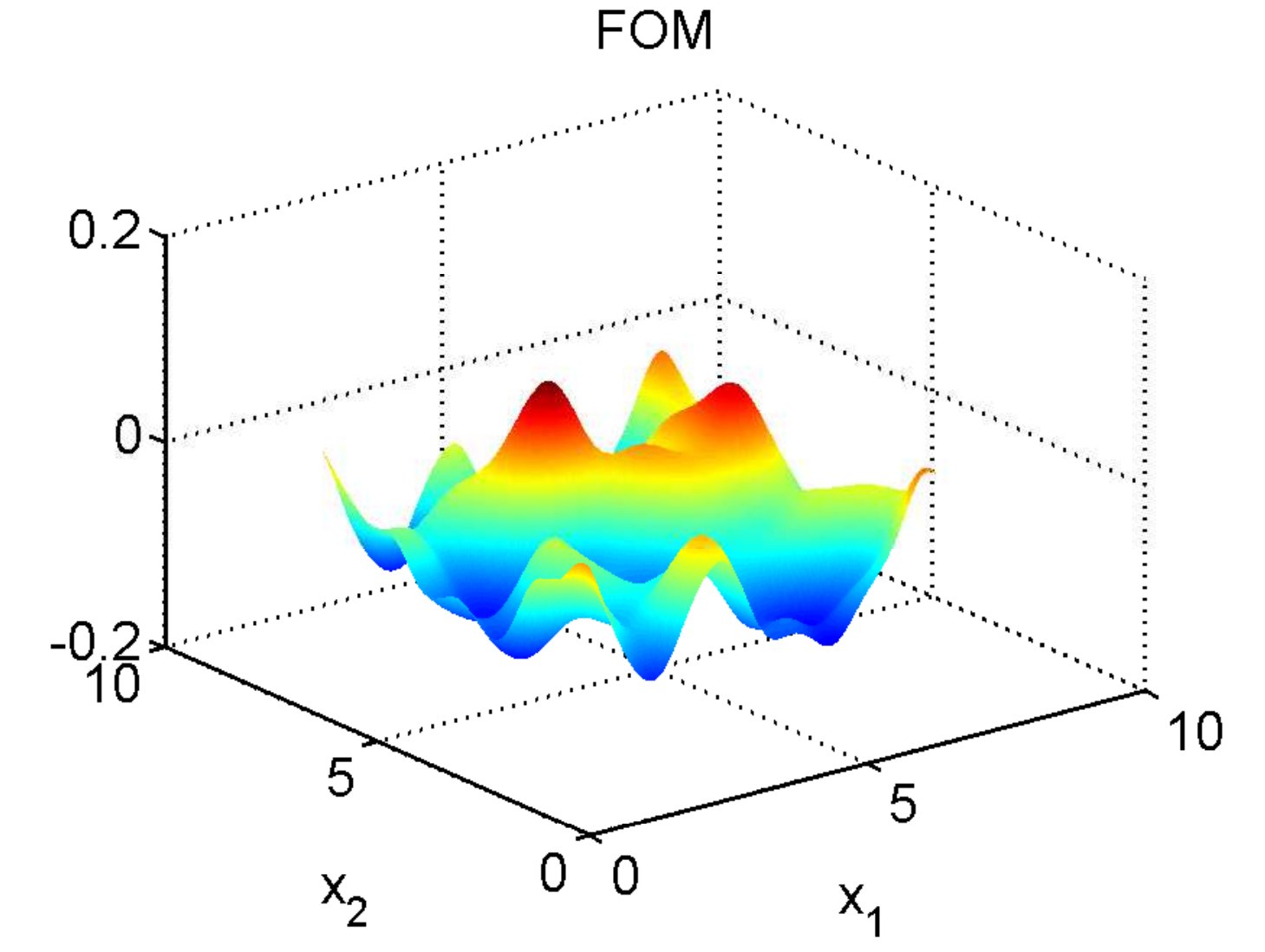}}
\subfloat{\includegraphics[width=0.45\textwidth]{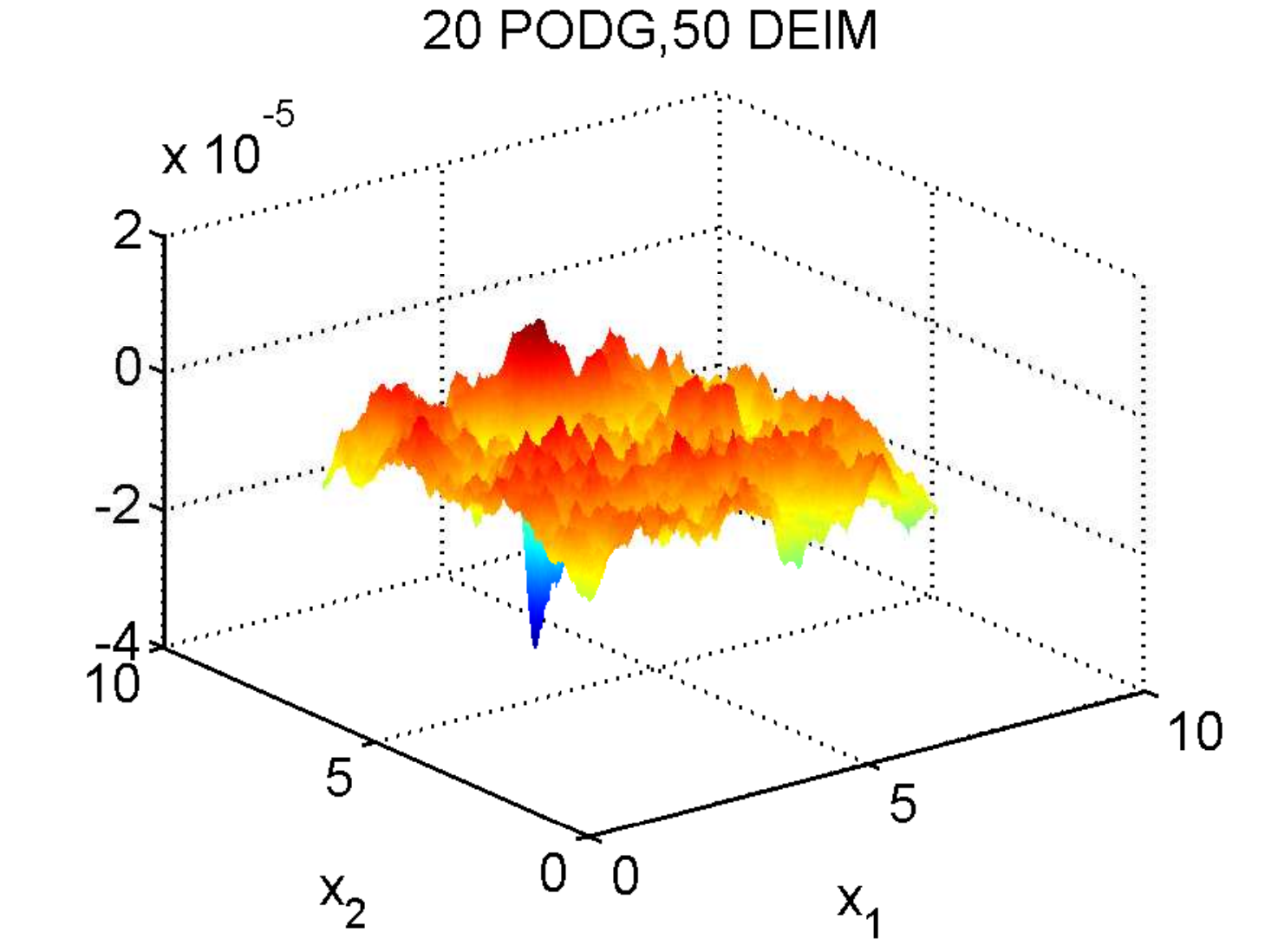}}
\caption{Allen-Cahn with logarithmic potential: FOM profile (left) and error plot between FOM and PODG-DEIM (right) solutions, at the final time for $\theta =0.10$.\label{ac_log_plot}}
\end{figure}

In the PODG algorithm, we choose now the temperature as a parameter by setting $\mu = \theta$ and we fix $\epsilon = 0.04$. The training set ${\mathcal M}_{train}$ is taken as the set consisting of the elements $\mu_k=0.05+0.03(k-1)\subset [0.05,0.17]$, $k=1,\ldots , 5$. The decrease of the error indicator is given in Fig.~\ref{ac_log_error}, left. The solution profile and error plot between FOM and PODG-DEIM solutions at the final time are presented in Fig.~\ref{ac_log_plot}, using 20 PODG and 50 DEIM modes. In Table~\ref{ac_table}, we give the numerical errors for the solutions and energies between FOM and ROM solutions for both examples.

\begin{table}[htb!]
\centering
\caption{FOM-ROM  solution errors and errors between discrete energies}
\label{ac_table}
\begin{tabular}{ l l c c c c }
  &  & \multicolumn{2}{c}{Quartic $F$} & \multicolumn{2}{c}{Logarithmic $F$}\\
  &  & Solution & Energy & Solution & Energy \\
\hline
PODG        &   & 9.87e-05 & 1.63e-06  & 9.66e-06  &  7.72e-07 \\
PODG-DEIM   &   & 9.94e-05 & 1.64e-06  & 1.08e-05  &  2.43e-06
\end{tabular}
\end{table}

\begin{table}[htb!]
\centering
\caption{CPU times and speed-up factors}
\label{speedup}
\begin{tabular}{ l c c c c c }
  &  \multicolumn{3}{c}{Wall Clock Time (sec)} & \multicolumn{2}{c}{Speed-Up Factor}\\
  & FOM & PODG & PODG-DEIM & PODG & PODG-DEIM \\
\hline
AC (quartic $F$)     &  50.64  & 9.72  & 2.23  &  5.21  &  22.71  \\
AC (logarithmic $F$) &  40.49  & 7.33  & 1.78  &  5.52  &  22.75
\end{tabular}
\end{table}

Finally, in Table~\ref{speedup}, we present the CPU times and speed-up factors related to the solutions of FOM, PODG without DEIM and PODG-DEIM. It can be easily seen that PODG with DEIM dramatically improves the computational efficiency.\\

\noindent {\bf Acknowledgments:}
The authors would like to thank the reviewer for the comments and suggestions that
helped to  improve the manuscript.



\end{document}